
\baselineskip=14pt
\parskip=10pt
\def\halmos{\hbox{\vrule height0.15cm width0.01cm\vbox{\hrule height
  0.01cm width0.2cm \vskip0.15cm \hrule height 0.01cm width0.2cm}\vrule
  height0.15cm width 0.01cm}}
\font\eightrm=cmr8 

\magnification=\magstephalf

\def\1{{\overline{1}}}
\def\2{{\overline{2}}}
\parindent=0pt
\overfullrule=0in

\def\frac#1#2{{#1 \over #2}}
\centerline
{\bf 
Counting Permutations Where The Difference Between Entries Located r Places Apart 
}
\centerline
{\bf
Can never be s (For any given positive integers r and s)
}
\bigskip
\centerline
{\it George SPAHN and Doron ZEILBERGER}

\bigskip

\qquad \qquad {\it In fond memory of David Peter Robbins (12 August 1942 - 4 September 2003) who was a great problem solver, and an equally good problem poser}

{\bf Abstract}: Given positive integers $r$ and $s$,  we use {\it inclusion-exclusion}, {\it weighted-counting of tilings}, 
and {\it dynamical programming}, in order to enumerate, semi-efficiently, 
the classes of permutations mentioned in the title.
In the process we revisit beautiful previous work of Enrique Navarrete, Roberto Tauraso, David Robbins (to whose memory this article is dedicated),
and John Riordan.
We also present two new proofs of John Riordan's recurrence (from 1965) for the  sequence enumerating  permutations without rising and falling
successions (the $r=1,s=1$ case of the title in the sense of absolute value).
The first is {\it fully automatic} using the (continuous) Almkvist-Zeilberger algorithm, while the second is {\it purely human-generated} via
an elegant combinatorial argument. We then state some open questions and pledge donations to the OEIS in honor of the solvers.
We conclude with a postscript describing a clever bijection by Rintaro Matsuo, that we were made aware of after the
first version was written, that suggests an alternative approach to our counting
problem, and that implies that there always exists, for each specific $r$ and $s$, a polynomial time algorithm for computing these sequences,
and also gives strong evidence (but not yet proves) some of the open problems. But as $r$ and $s$ get larger, the implied degree of the
``polynomial time algorithm" grows large enough to guarantee that our original, ``exponential time'' algorithms, are still useful.

{\bf Important: Maple Packages and Output Files}

This article is accompanied by two Maple packages, {\tt ResPerms.txt}, and {\tt ResPermsRM.txt}, downloadable, respectively, from

{\tt https://sites.math.rutgers.edu/\~{}zeilberg/tokhniot/ResPerms.txt} \quad , and

{\tt https://sites.math.rutgers.edu/\~{}zeilberg/tokhniot/ResPermsRM.txt} \quad .

The front of this article 

{\tt https://sites.math.rutgers.edu/\~{}zeilberg/mamarim/mamarimhtml/perms.html} \quad,

contains many input and output files.

{\bf How it all Started}

A few months ago we came across a lovely paper by Manuel Kauers and Christoph Koutschan [KK] where they described
how to guess recurrences satisfied by `hard to compute sequences', where one gets very plausible conjectured
recurrences for sequences with far fewer ``data points" than are needed by traditional guessing with
{\it undetermined coefficients}. As a very impressive {\it case study} they found a {\it conjectured} (but {\it practically certain})
{\bf linear recurrence with polynomial coefficients} for the sequence enumerating 
permutations $\pi$ of $\{1, \dots, n\}$ such that for all $1 \leq i \leq n-2$, one has:
$$
\pi_{i+2} - \pi_i \neq 2 \quad .
$$

This is sequence  A189281 in the OEIS ({\tt https://oeis.org/A189281}) [S22].
Initially the OEIS only had $35$ terms, computed with great computing effort by  Vaclav Kotesovec,  hence it was labeled {\bf hard}.
There was no {\it a priori} theoretical guarantee that it is {\it holonomic}, i.e. satisfies a linear recurrence equation with polynomial coefficients, 
and for that matter, we still don't have one, but {\it surprisingly}, even with this `little data', Kauers and Koutschan, using their
new method, were able to conjecture a linear recurrence equation of order $8$ and degree $11$.
In general, using the traditional {\it vanilla} method of {\it undetermined coefficients},
in order to conjecture a recurrence of order $r$ and degree $d$ requires at least $(r+1)(d+1)+r+1$ data points,
hence to conjecture a recurrence of that order and degree would require $(8+1)\cdot (11+1)+8+1=117$  data points, rather than the at-the-time available $35$ terms.
Nevertheless, they succeeded using so little data!

Assuming that the [KK] conjecture was true, one of us (DZ) found an equivalent recurrence of order $13$ and degree $3$ that would only necessitate $(13+1)\cdot(3+1)+13+1=70$ data points (see [S22]).

But even with the new method of [KK], one needs efficient ways to generate as many terms as possible. 
Our goal is to compute many  terms of the following two doubly-infinite families of sequences,
for {\it any given} positive integers $r$ and $s$
$$
a_{r,s}(n):= \# \{\pi \in S_n \, ; \quad \pi_{i+r}\,- \, \pi_r \neq \, s \,\, for \quad all \quad 1 \leq i \leq n-r\} \quad,
$$
and
$$
b_{r,s}(n):= \#\{\pi \in S_n \, ; \quad |\pi_{i+r}\,- \, \pi_r| \neq \, s \,\, for \quad all \quad 1 \leq i \leq n-r\} \quad .
$$

(For any finite set $S$, $\#S$ denotes its number of elements).

Using our method, we were able to compute {\it from scratch}, $70$ terms of {\tt A189281}, so in {\it hindsight}, one did not need the new Kauers-Koutschan method for this particular
sequence, but we are sure that it would be useful in many other cases.

{\bf Past Work}

For the sequences $a_{1,s}(n)$, for {\it any positive integer} $s$, Enrique Navarrete found a beautiful and simple explicit formula, using a clever {\it inclusion-exclusion} argument.

{\bf Theorem} (E. Navarrete [N]): for all $s \geq 1$ and $n \geq s$,
$$
a_{1,s}(n) \, = \,\sum_{j=0}^{n-s} (-1)^j {{n-s} \choose {j}} (n-j)! \quad .
$$

Using the {\bf Zeilberger algorithm}[Z1] we immediately get the following corollary.

{\bf Corollary:} For $s \geq 1$ and for $n \geq 2$ we have,
$$
a_{1,s}(n)=(n-1)\, a_{1,s}(n-1)+(n-s-1)\,a_{1,s}(n-2) \quad .
$$

Later we will give a combinatorial proof of this recurrence for $a_{1,1}(n)$.

The sequence $b_{1,1}(n)$ is famous! It is sequence A2464 [S11a] {\tt https://oeis.org/A002464}, called         
Hertzsprung's problem, and was already investigated by Paul Poulet [P] in 1919. Its description in the OEIS is as follows:

{\it ``Ways to arrange $n$ non-attacking kings on an $n \times n$ board, with $1$ in each row and column. Also number of permutations of length $n$ without rising or falling successions.''}

The combinatorial giant, John Riordan [Rio],  proved a nice fourth-order recurrence
$$
b_{1,1}(n)\,= \,(n+1) \, b_{1,1}(n-1) \, - \, (n-2)\, b_{1,1}(n-2) \, - \, (n-5)\,b_{1,1}(n-3) \, + \, (n-3)\,b_{1,1}(n-4) \quad .
$$
We will later give two new proofs. The first using the continuous Almkvist-Zeilberger algorithm [AZ], and the second using an elegant combinatorial argument.

Our hero, Dave Robbins [Ro], used  a clever {\it inclusion-exclusion} argument to prove the following double sum formula for $b_{1,1}(n)$:
$$
b_{1,1}(n)=  \sum_{i=0}^{n-1} (-1)^i (n-i)! \, \sum_{c=1}^{i} {{i-1} \choose {i-c}} {{n-i} \choose {c}}2^c \quad .
$$

We will later use Dave Robbins' approach to study, and efficiently compute, the sequences $b_{1,s}(n)$ for {\it any} $s>1$, none of which are (yet) in the OEIS.

To conclude this historical section let's mention which sequences are currently (Dec. 2022) in the OEIS:

$\bullet$ $a_{1s}(n)$ for $1 \leq s \leq 5$, see references [S11],[S12],[S13],[S14],[S15] respectively.

$\bullet$ $a_{rr}(n)$ for $1 \leq r \leq 6$, see references [S11],[S22],[S33],[S44],[S55],[S66] respectively.

$\bullet$ $b_{rr}(n)$ for $1 \leq r \leq 6$, see references [S11a],[S22a],[S33a],[S44a],[S55a],[S66a] respectively.

None of the  sequences $a_{rs}(n)$  and  $b_{rs}(n)$  with $r >1, s>1$ and $r \neq s$ are yet (Oct. 2022) in the OEIS.
Neither are the sequences $b_{1,s}(n)$ for $s>1$, for which we will later give polynomial-time algorithms.

In the  next section we will show how to compute many terms of these new sequences and compute {\it many more terms} for the
already present sequences $a_{rr}(n)$ and $b_{rr}(n)$, $2\leq r\leq 6$.

{\bf Computing as Many Terms as Possible of the sequences ${\bf a_{r,s}(n)}$ and ${\bf b_{r,s}(n)}$ }

We will use the old workhorse of {\it inclusion-exclusion} but with a new {\it twist}, that would make it amenable for {\it symbolic computation}.

Fix $r \geq 1$ and $s \geq 1$. We want to count all the {\it good guys}, i.e. permutations of $\{1, \dots , n\}$ such that none of the following $n-r$ {\bf unfortunate events} are
committed
$$
\pi_{r+1} -\pi_1=s \quad, \quad \pi_{r+2} -\pi_2=s \quad, \quad \dots \quad, \quad \pi_{n}-\pi_{n-r}=s \quad .
$$
As usual, instead of counting {\bf good guys}, we do a {\it signed} counting of {\bf all pairs}, 
$$
[guy, S] \quad,
$$
where {\it guy} is {\it any} permutation, and $S$ is a subset (possibly empty, possibly the whole) of its set of unfortunate events. (See [Z3] for an engaging account).
Each such pair contributes $(-1)^{\#S}$ to the total sum. While this new signed sum has many more terms, and is extremely inefficient when applied to
specific sets, it is a great {\bf theoretical tool} when used cleverly by humans (and computers!).

Looking at the structure of such a {\it set of unfortunate events}, we see that the {\bf board} $\{1,2, \dots, n\}$ has a certain subsest of the set of pairs
$$
\{\{1,r+1\} \, , \, \{2,r+2\} \, , \quad \dots , \, \{n-r,n\} \}\quad,
$$
Think of them as a collection of {\it arcs} whose width is $r$ . Each such configuration naturally forms into a disjoint
union of {\bf connected components}. All the entries that do not participate form singleton tiles. So looking at the possible connected components,
we see a natural {\it tiling} of the `board'  $\{1,2, \dots, n\}$ into  horizontal shifts of tiles of the form
$$
\{1\} \, , \, \{1, 1+r \} \, , \, \{1, 1+r , 1+ 2r\} \, , \, \{1, 1+r , 1+ 2r, 1+3r\} \, , \,  \dots \quad .
$$

So let's introduce variables $x_1 , x_2 , x_3 \dots $ and declare that the weight of a tile $t$ is $x_{\#t}$. In particular, the weight of a singleton is $x_1$.
The weight of a {\it tiling} is the {\bf product} of the weights of all its constituent tiles. Let's give a simple example.

Suppose that $r=2$  and $s=3$,  $n=9$ and the subset $S$, of  the set of $\pi$'s unfortunate events, happens to be:
$$
\{\pi_{3}-\pi_1=3 \quad, \quad \pi_{5}-\pi_3=3 \quad, \quad \pi_9-\pi_7=3 \} \quad .
$$
Since $\{1,3\}$ and $\{3,5\}$  form a {\bf connected component} $\{1,3,5\}$, and $\{7,9\}$ forms its own connected component, and
the integers not participating in these two connected components are all singletons, we have the tiling
$$
\{ \{1,3,5\} \, , \, \{7,9\} \, , \, \{2\} \, , \, \{4\}  \, ,  \, \{6\} \, , \, \{8\} \} \quad .
$$

So let's define a polynomial in the (potentially infinite, but for any given $n$ finite) set of variables $x_1,x_2, x_3, \dots$
$$
f_{r,n}(x_1, x_2, x_3 \dots ) \, := \, \sum_{T} weight(T) \quad,
$$
where the sum is taken over all the tilings of the board $\{1, 2, \dots, n\}$ by the tiles (i.e. horizontal shifts) of
$$
\{1\} \quad , \quad \{1,1+r\} \quad , \quad \{1,1+r, 1+2r\} , \quad \{1,1+r, 1+2r, 1+3r\} \quad , \quad \dots.
$$

For example
$$
f_{3,5}(x_1,x_2,x_3)\, = \, x_{1}^{5}+2 x_{1}^{3} x_{2}+x_{1} x_{2}^{2} \quad,
$$
since 

$\bullet$ There is one tiling with all singletons (horizontal shifts of $\{1\}$): 
$$
\{\{1\}  , \{2\} , \{3\}  , \{4\}  , \{5\}  \} \quad,
$$
whose weight is ${x_1}^{5}$. 

$\bullet$ There are two tilings with $3$ singletons and one horizontal shifts of $\{1,4\}$, namely
$$
\{ \{1,4\}, \{2\} , \{3\}  , \{5\}  \}  \quad, \quad
\{  \{2,5\}, \{1\},  \{3\}, \{4\}  \}  \quad,
$$
each of whose weight is ${x_1}^{3}x_2$. 

$\bullet$ There is one tiling with $1$ singleton and two  horizontal shifts of $\{1,4\}$, namely
$$
\{ \{1,4\},  \{2,5\} , \{3\} \}\quad ,
$$
whose weight is $x_1 {x_2}^{2}$. 

Using {\it dynamical programming} (see [Z2]), it is very fast to compute these polynomials. Since these {\it weight-enumerators} of tilings are so {\it central} to our approach, we
invite our readers to convince themselves that
$$
f_{3,7}(x_1,x_2,x_3)\, = \,
x_{1}^{7}+4 x_{1}^{5} x_{2}+x_{1}^{4} x_{3}+5 x_{1}^{3} x_{2}^{2}+2 x_{1}^{2} x_{2} x_{3}+2 x_{1} x_{2}^{3}+x_{2}^{2} x_{3} \quad.
$$

{\eightrm This is implemented in procedure {\tt fsn(r,n,x)} in our Maple package.}

But there is {\bf another tiling} present. If you look at the {\bf actual values} of the `unfortunate events', you have a natural tilings with tiles
$$
\{1\} \, ,\, \{1, 1+s \} \, , \, \{1, 1+s , 1+ 2s\} \, , \, \{1, 1+s , 1+ 2s, 1+3s\}\, , \dots  \quad,
$$

whose weight-enumerator is $f_{s,n}(x_1, x_2, \dots )$.

For example, look at one of the many  pairs that come up in trying to compute $a_{2,3}(9)$:
$$
[326195487, \{ \pi_3-\pi_1=3, \pi_5-\pi_3=3, \pi_9-\pi_7=3 \}] \quad .
$$
Looking at the {\bf values} we have the tiling (since $\pi_1=3,\pi_3=6,\pi_5=9$ and $\pi_7=4, \pi_9=7$):
$$
\{ \, \{3,6,9\} \, , \,  \{4,7\} \, ,  \, \{1\},\{2\},\{5\},\{8\} \, \}  \quad ,
$$
whose weight is the same as  that of the corresponding {\it input tiling}, $\{ \, \{1,3,5\} \, , \,  \{7,9\} \, ,  \, \{2\},\{4\},\{5\},\{8\} \, \}$,
namely  $x_1^4 x_2 x_3$.

Now write
$$
f_{r,n}(x_1, \dots , x_n) \, = \, \sum C^{(n,r)}_{\alpha} x_1^{a_1} x_2^{a_2} \cdots x_n^{a_n} \quad,
$$
where the sum ranges over all {\bf integer partitions} $\alpha=1^{a_1} 2^{a_2} \dots n^{a_n}$ of $n$, written in {\bf frequency notation}. Note that, of course: $1\cdot a_1 + 2 \cdot a_2 + \dots + n\cdot a_n=n$.
(Also note that  if $r=1$ all integer partitions show up, but when $r>1$ then some never do).

{\bf Finally} we are ready to state our `formula' for $a_{r,s}(n)$.
$$
a_{r,s}(n)= \sum_{\alpha}\,  C^{(n,s)}_{\alpha} \cdot C^{(n,r)}_{\alpha} \cdot (-1)^{ a_1+a_2+ \dots +a_n -n} a_1! a_2! \cdots a_n! \quad,
$$
where the sum ranges over all partitions of $n$ $\alpha=1^{a_1} 2^{a_2} \dots n^{a_n}$, written in frequency notation. 

Let's explain. Any set of unfortunate events gives rise to {\bf two} tilings: 

$\bullet$ One with tiles $\{1\}, \{1,1+r\}. \{1,1+r,1+2r\}, \dots$ taken care by the $f_{r,n}$ ; and

$\bullet$ one with tiles $\{1\}, \{1,1+s\}. \{1,1+r,1+2s\}, \dots$  taken care by the $f_{s,n}$ . 

The cardinality of the set of resulting unfortunate events is
$$
a_1 \cdot (1-1)+ a_2 \cdot (2-1)+ a_3 \cdot (3-1) + \dots = n - (a_1+ \dots a_n) \quad,
$$
explaining the exponent of $(-1)$, that comes from inclusion-exclusion. 
Finally, one has to match the {\it `input tiles'} and the {\it `output tiles'}. Of course they have to be of compatible sizes.

$\bullet$  For each $i$, there are $a_i!$ ways to match the $a_i$  tiles of size $i$ of the input tiling with the $a_i$ tiles of the output tiling.

The total number of possible {\it matches} is thus:
$$
a_1! a_2! \cdots a_n! \quad .
$$

So we explained all the ingredients in our `formula'.

A very small tweak gives the formula for $b_{r,s}(n)$:
$$
b_{r,s}(n)= \sum_{\alpha}\, C^{(n,s)}_{\alpha} \cdot C^{(n,r)}_{\alpha} \cdot (-1)^{ a_1+a_2+ \dots +a_n -n} a_1! a_2! \cdots a_n! \cdot 2^{a_2+a_3+ \dots +a_n} \quad,
$$
where the sum ranges over all partitions of $n$, $1^{a_1} 2^{a_2} \dots n^{a_n}$, written in frequency notation.

The extra factor of  $2^{a_2+a_3+ \dots +a_n}$ in the summand comes from the fact that now we 
count permutations of $\{1, \dots, n\}$ such that $|\pi_{i+r}-\pi_i| \neq s$, hence now in each component of the tiling that is not a singleton,
there are {\bf two} possible directions, {\bf up} and  {\bf down}.

{\eightrm The formulas for $a_{r,s}(n)$ and  $b_{r,s}(n)$  are implemented by procedures {\tt Urs(r,s,n)} and {\tt Vrs(r,s,n)}, respectively, in our Maple package {\tt ResPerms.txt} .}

{\bf Did we answer the enumeration question?}

Well, not quite! According to Herb Wilf's classic essay [W] we need a {\bf polynomial time algorithm} and the number of terms in our `formulas' for $a_{r,s}(n)$ and $b_{r,s}(n)$ is the
number of partitions of $n$ (in fact slightly less, but not significantly less, asymptotically speaking). According to Hardy and Ramanujan they
are about $e^{\pi\sqrt{2n/3}}$ of them. Nevertheless our approach worked very well up to $n=60$ and enabled us to get, with relatively little computing effort, many more terms.
These sequences were first computed by Vaclav Kotesovec, that according to his OEIS {\it profile}, is
a FIDE International master in Chess Problems. (Kotesovec gives nice interpretations of these sequences in terms of  ``non-attacking" placements of some {\it fairy-chess} pieces).

Just to illustrate the (relative) efficiency of our new method,
in order to get the first $30$ terms of the sequence $a_{4,4}(n)$ ({\tt https://oeis.org/A189283}), we only needed {\bf less than three seconds}, on our modest laptop. Here there are:
$$
1, 2, 6, 24, 114, 628, 4062, 30360, 255186, 2414292, 25350954, 292378968, 3673917102, 49928069188, 
$$ 
$$
729534877758, 11403682481112, 189862332575658, 3354017704180052, 62654508729565554, 1233924707891272728, 
$$
$$
25550498290562247438, 554913370184289495780, 12612648556263898345758, 299411750583810718488216, 
$$
$$
7409924986737790240296258, 190856850583975937020030228, 5108283222440036893650974970, 
$$
$$
141870112250977140975169694808, 4082973503947066134710463043374, 121616802487841972048586204012740 \quad .
$$

For many terms of such sequences, including $a_{r,s}(n)$ and $b_{r,s}(n)$ with $r>1, s>1$ and $r \neq s$, none of which are yet in the OEIS, see the numerous output files in the front of this article:

{\tt https://sites.math.rutgers.edu/\~{}zeilberg/mamarim/mamarimhtml/perms.html}.

{\bf Comment}: After writing the first version, there were some interesting developments made by Rintaro Matsuo [Rin].
See the postscript.

{\bf The beautiful work of Roberto Tauraso}

We should mention that the case of $b_{r,r}(n)$ for $r=2$ and (larger $r$) has been nicely handled by Roberto Tauraso [T]. We believe that his
approach bears some similarities to ours, but ours is more amenable to using symbolic computations efficiently.

In particular the sequence $b_{2,2}(n)$, currently only has $35$ terms listed in {\tt https://oeis.org/A110128}. See our output file

{\tt https://sites.math.rutgers.edu/\~{}zeilberg/tokhniot/oResPerms22Va.txt} 

for $60$ terms!

{\bf Update}: For $150$ terms, see the Postscript.

{\bf Polynomial-Time Algorithms for computing $b_{1,s}(n)$ for all $s$}

As pointed out by Enrique Navarrete, the formula for $a_{1,s}(n)$  is very simple, and as we already mentioned, it implies (via the Zeilberger algorithm [Z1]), a very simple
{\bf second-order} recurrence. Dave Robbins [Ro] (and before him John Riordan [Rio] and even before him, Paul Poulet [P]) explored the sequence $b_{1,1}(n)$ also using
{\bf inclusion-exclusion}, but this time there are no more (potentially) {\it infinitely-many} tiles to keep track of, 
because the tiles in the `output' domain are very simple, they are just sets of the form $\{1\}, \{1,2\}, \{1,2,3\}, \dots$, so
one only needs to find the weight-enumerator of tiles according to the weight $x^{TotalNumberOfTiles}$, without recording
the individuality of the tiles. {\bf But} one has to keep track of the {\it total length} of the non-singleton tiles,
so we need to introduce another formal variable, let's call it $z$ that keeps track of the total length of the non-singleton tiles. 
Also  each non-singleton tile can have {\bf two directions} so we would have to replace $z$ by $2z$.
The details may be gleaned from the {\it Maple source code} of procedures {\tt Ker1s(s, x, z, K)}.
We believe that it can be proved that these Robbins-style weight-enumerators of tilings are
{\bf always}  {\it rational functions}, but being {\it experimental mathematicians}, we just guessed them, so {\it officially}
the output of the recurrences for the sequences for $b_{1,s}(n)$ with $s>1$ are only {\it semi-rigorous}. 
(On the other hand, thanks to Dave Robbins, our automatic proof of Riordan's recurrence for the case $s=1$, i.e. for $b_{1,1}(n)$, is {\it fully rigorous}.)

{\it At  the end of the day}, once you have the Robbins-weight-enumerator $K_s(x,z)$ for enumerating $b_{1,s}(n)$, that is a rational function, one
applies the {\it umbra} $z^n \rightarrow n!$, but this is {\it the same as} multiplying by $e^{-z}$,  and then integrating from $0$ to $\infty$.

Hence we have the following theorem.

{\bf Semi-Rigorous Theorem:} For $s \geq 1$, there almost exists a rational function of $x$ and $z$ (that Maple can find), let's call it $K_{s}(x,z)$, such that
$$
\sum_{n=0}^{\infty} \, b_{1,s}(n) x^n \, = \, \int_0^{\infty} K_s(x,z) e^{-z} \, dz \quad .
$$
Using the {\bf Almkvist-Zeilberger} algorithm [AZ] one can find an (inhomogeneous)  {\bf linear differential equation} with {\bf polynomial coefficients},
that translates to a (homogeneous) {\bf linear recurrence with polynomial coefficients} for the actual sequence $b_{1,s}(n)$.

See the output files

{\tt https://sites.math.rutgers.edu/\~{}zeilberg/tokhniot/oResPermsR1.txt} \quad,

[This output file contains the first promised new proof of Riordan's recurrence, using the Almkvist-Zeilberger algorithm].

{\tt https://sites.math.rutgers.edu/\~{}zeilberg/tokhniot/oResPermsR2.txt} \quad,

{\tt https://sites.math.rutgers.edu/\~{}zeilberg/tokhniot/oResPermsR3.txt} \quad,

{\tt https://sites.math.rutgers.edu/\~{}zeilberg/tokhniot/oResPermsR4.txt} \quad,

for the sequence $b_{1,1}(n)$, $b_{1,2}(n)$, $b_{1,3}(n)$, $b_{1,4}(n)$, respectively.

{\bf  Combinatorial, Human-Generated, Proofs of the Linear Recurrences satisfied by $a_{1,1}(n)$ and $b_{1,1}(n)$}.

{\bf Proposition 1}: Let $a(n) = a_{1,1}(n)$. Then $$ a(n) = (n-1) \cdot a(n-1)+(n-2) \cdot a(n-2) \quad .$$

{\bf Combinatorial Proof}: We would like to count the number of permutations $\pi$ that never have $\pi_{i+1}-\pi_i=1$. Suppose we had such a permutation $\alpha$, and we remove the $n$. We now have a permutation $\beta$ of length $n-1$. There are two cases. Either $\beta$ has no $i$ such that $\beta_{i+1}-\beta_i=1$, or it has a violation, an $i$ such that $\beta_{i+1}-\beta_i=1$. Since $\alpha$ had no violations, the only possible violation in $\beta$ is involving the 2 elements that used to be separated by $n$.

{\bf Case 1}: There is no violation in $\beta$. There are $a(n-1)$ possibilities for $\beta$ in this case. In each one, there are $n$ possible locations we could insert $n$ to get a permutation of length $n$. As long as we don't insert the $n$ immediately after the $n-1$, we will get a permutation of length $n$ with no violations. Thus for each possible $\beta$, there are $n-1$ possible $\alpha$ that could map to it. Thus the number of permutations in this case is equal to $(n-1) \cdot a(n-1)$.

{\bf Case 2}: Removing $n$ from $\alpha$ caused a single violation in $\beta$. Let us define an auxiliary sequence, $b(n)$, which counts the number of permutations $\pi$ of length $n$ with a single $i$ such that $\pi_{i+1}-\pi_i=1$. Each possible $\pi$ in this case, can be made into a permutation of length $n$ with no violations by inserting $n$ in between the 2 elements involved in the violation. Thus $b(n-1)$ counts the number of possibilities for $\alpha$ in case 2. 

We now have that $$ a(n) = (n-1) \cdot a(n-1)+b(n-1) \, .$$
We next argue that $$b(n) = (n-1) \cdot a(n-1) \, .$$ 
Consider the operation $f_i$ for $1 \leq i \leq n$ that maps permutations $\pi$ of length $n$ to permutations of length $n+1$. $f_i(\pi)$ is computed in the following way. First insert an $i+1$ in $\pi$ immediately after the $i$. Then relabel all elements $j > i$ as $j+1$. For example, $f_2(321) = 4231$. If we apply one of these operations $f_i$ to a permutation with no violations, we get a permutation with exactly one violation. Further, any permutation with exactly one violation can be obtained by applying an $f_i$ to a permutation with no violations. Thus the set of permutations described by $b(n)$ are in bijection with $\{f_i(\pi)\}$, where $\pi$ a permutation described by $a(n-1)$. There are $n-1$ possibilities for $i$ and $a(n-1)$ possibilities for $\pi$, so we get that indeed $b_n = (n-1) \cdot a(n-1)$. Thus $$b(n-1) = (n-2) \cdot a(n-2) \, ,$$ and we get $$ a(n) = (n-1) \cdot a(n-1)+(n-2) \cdot a(n-2)\, . \quad \halmos$$ 

In a  talk ({\tt https://www.youtube.com/watch?v=ahuY7FHPiAc}) delivered at
the conference  ICECA 2022 (September 6-7, 2022) ({\tt http://ecajournal.haifa.ac.il/Conf/ICECA2022.html}), the second-named author (DZ) pledged to donate $\$100$ to the OEIS in honor of
the first person to find a {\bf combinatorial} proof of Proposition 2 below. This challenge was met by the first-named author (GS). The donation was made.

{\bf Proposition 2}: Now let $a(n) = b_{1,1}(n)$. Then $$ a(n) = (n+1) \cdot a(n-1)-(n-2) \cdot a(n-2)-(n-5) \cdot a(n-3)+(n-3) \cdot a(n-4) $$

{\bf Combinatorial Proof}: In this case we say a up-violation in a permutation $\pi$ is an integer $i$ so that $\pi_{i+1}-\pi_i=1$ and a down-violation is an integer $i$ so that $\pi_{i+1}-\pi_i=-1$. We now define 2 auxiliary sequences, $b(n)$ and $c(n)$.

$A(n)$ is the set of permutations of length $n$ with no violations (of either kind). 

$a(n) = \# A(n)$

$B(n)$ is the set of permutations of length $n$ with exactly one violation (of either kind). 

$b(n) = \# B(n)$

$C(n)$ is the set of permutations of length $n$ with exactly one violation (of either kind) at specifically $i= n-1$.

$c(n) = \# C(n)$

{\bf Lemma 1}: $$ a(n) = (n-2) \cdot a(n-1) +b(n-1) -c(n-1) \, .$$

For the proof, we first state that if you start with a permutation of length $n-1$ with no violations, then there are $n-2$ locations that we can insert $n$ into such that we still do not have any violations. Doing so only requires that we avoid inserting $n$ next to $n-1$. We can generate $(n-2) \cdot a(n-1)$ unique permutations in this manner. This will generate some of the permutations in $A(n)$ but not all. The others are permutations $\pi$ such that if $n$ is removed, we introduce a single violation. Removing $n$ from $\pi$ will give an element of $B(n-1)$, but not all elements of $B(n-1)$ can be reached. If we have an element of $B(n-1)$ that is also an element of $C(n-1)$, then inserting $n$ will introduce a new violation, so we will not get a permutation in $A(n)$. In summary, any permutation in $B(n) \setminus  C(n)$ can be expanded into an element of $A(n)$ by inserting $n$ between the violation. This gives all remaining elements of $A(n)$, so we get that  $ a(n) = (n-2) \cdot a(n-1) +b(n-1) -c(n-1) $

{\bf Lemma 2}: $$ b(n) = 2(n-1) \cdot a(n-1) +2b(n-1) +b(n-2) $$

We can obtain most of the elements of $B(n)$ in the following way. Start with an element of $A(n-1)$, $\pi$. As in the proof of Proposition 1, we will use an operation similar $f_i$ to expand one elements $i$ into either $i, i+1$ or $i+1,i$, and increase the value of each element $j \in \pi$, $j > i$ by 1. So the function outputs expanding $i$ into an up-violation or a down-violation. For example, $f_1(53142) =\{641253 , 642153\}$. Then there $(n-1)$ choices for $i$, $a(n-1)$ choices for $\pi$, and 2 choices for whether it expands into an up-violation or down-violation. This gives $2(n-1) \cdot a(n-1)$ elements of $B(n)$, and like last time it only remains to count the ones where the inverse of this procedure introduces new violations. For example $4213$ and $4132$  would collapse to $312$ and not be generated by an element of $A(n-1)$. These elements collapse to a permutation with at least one violation. The ones that collapse to having exactly one violation end up in $B(n-1)$. Each permutation $\alpha$ in $B(n-1)$ can be expanded into a permutation in $B(n)$ in $2$ ways. Suppose $\alpha$ has an up-violation at $i$. Then we replace $i$ with the down-violation $i+1, i$, and increase all $j > i$ by 1. This necessarily fixes the previous violation, and does not introduce any violations. Alternatively, we could replace $i+1$ with $i+2,i+1$, which would have the same effect. If instead $\alpha$ had a down violation, $i+1,i$, you could replace it with $i+2,i,i+1$, or $i+1,i+2,i$. There are always exactly 2 permutations in $B(n)$ which map down to each permutation in $B(n-1)$, so we get $2b(n-1)$ contributions from this case. We also have to account for the possibility that after collapsing we have both new possible violations. For example $1324$ would collapse to $123$ which is not in $B(n-1)$. Luckily we can further collapse this triple violation and end up with a permutation in $B(n-2)$. Thus $ b(n) = 2(n-1) \cdot a(n-1) +2b(n-1) +b(n-2) $.

{\bf Lemma 3}: $$ c(n) = 2 a(n-1) +c(n-1) \, .$$

For this one, we can start with an element of $A(n-1)$ and insert the $n$ next to the $n-1$ to get an element of $C(n)$. There are 2 places we can insert it, so we get $2 \cdot a(n-1)$  unique elements of $C(n)$ in this way. We then have the elements of $C(n)$ where removing $n$ causes $n-1$ to be adjacent to $n-2$. The number of permutations in this case is exactly $c(n-1)$. All elements of $C(n)$ fall into one of the two cases, so we get $ c(n) = 2 a(n-1) +c(n-1) $. 

The above 3 lemmas give a system of 3 recurrences and 3 equations. We can solve this using ({\bf non-commutative}) linear algebra. Letting $N$ be the {\bf negative} shift operator, 
$$
N\, f(n) := f(n-1) \quad,
$$
and using $a,b,c$ to refer to the sequences as a whole, we get:

$$ a = (n-2) N a + N b - N c $$
$$ b = 2(n-1) N a + 2 N b + N^2 b $$
$$ c = 2 N a + N c$$

Solving for $a$ using Gaussian elimination, and the fact that  $N n=(n-1) N$, gives $$ [ 1-(n+1) +(n-2)N^2+(n-5)N^3-(n-3)N^4 ] a = 0 \, ,$$
so we conclude that
$$ a(n) = (n+1) \cdot a(n-1)-(n-2) \cdot a(n-2)-(n-5) \cdot a(n-3)+(n-3) \cdot a(n-4) \, . \quad \halmos $$

{\bf Final Challenges}

{\bf Challenge 1} ($\$ 100$ to the OEIS in honor of the first prover or disprover) Rigorously prove the
recurrence conjectured in [KK] for the sequence $a_{2,2}(n)$ (or the equivalent one found by DZ) given in [S22].

{\bf Challenge 2} ($\$ 200$  to the OEIS in honor of the first conjecturer)  Using the method of [KK] or otherwise,
{\it conjecture} a recurrence for the sequence $b_{2,2}(n)$, whose first $60$ terms are given here:

{\tt https://sites.math.rutgers.edu/\~{}zeilberg/tokhniot/oResPerms22Va.txt} \quad .

(Of course you are welcome to generate more terms! Update: see the postscript for $150$ terms). 

{\bf Update:} Rintaro Matsuo (see the postscript) used his method to generate $600$ terms! This enabled Manuel Kauers and Christoph Koutschan to meet this challenge.
A donation to the OEIS in honor of Matsuo, Kauers, and Koutschan has been made.

{\bf Challenge 2'} (Additional $ \$ 100$). Prove it rigorously.

{\bf Challenge 3} ($\$ 300$ to the OEIS in honor of the first prover or disprover) Prove (or disprove) that for every $r>1$ and $s>1$,
the sequences $a_{r,s}(n)$ and $b_{r,s}(n)$ are {\it always} {\it holonomic} (alias {\it P-recursive}), in other words are satisfied by some {\it homogeneous linear recurrence equation with polynomial coefficients}.
[Currently (as far as we know) there is no `general theory' that would imply this.]

\vfill\eject

{\bf Postscript: Rintaro Matsuo's Brilliant Bijection}

After the first version of this article was written, we were made aware (via the OEIS) that Rintaro Matsuo [Rin] found a {\it polynomial time} algorithm 
($O(n^4)$ for each single term, and hence $O(n^5)$ for the first $n$ terms),
to compute the sequence $a_{2,2}(n)$, that was the  {\it starting point} of our paper. He also implemented it and produced $300$ terms of
{\tt https://oeis.org/A189281}.

The key insight of Matsuo is the bijection of the set $\{1, \dots, n\}$ to itself given by
$$
(1,2,3, \dots, n)  \, \rightarrow \, (
1 \, , \,  1+ \lfloor \frac{n+1}{2} \rfloor \, , \,
2\, ,\, 2+ \lfloor \frac{n+1}{2} \rfloor \,,\,
3 \, , \, 3+ \lfloor \frac{n+1}{2} \rfloor, \dots )
$$

that transforms the set of permutations of length $n$ such that $\pi_{i+2}-\pi_i \neq 2$ to a {\it superset} of the much easier set of
permutations such that $\pi_{i+1}-\pi_i \neq 1$, except that one is allowing  $\pi_{i+1}-\pi_i =1$ if $i= \lfloor \frac{n+1}{2} \rfloor$
or $\pi_i=\lfloor \frac{n+1}{2} \rfloor$. The Maple package {\tt ResPermsRM.txt}, available from

{\tt https://sites.math.rutgers.edu/\~{}zeilberg/tokhniot/ResPermsRM.txt} \quad ,

gives a Maple implementation of Rintaro Matsuo's algorithm [procedure {\tt U22rm(n)}]. For the first $150$ terms see:

{\tt https://sites.math.rutgers.edu/\~{}zeilberg/tokhniot/oResPerms22rmA.txt} \quad .

However, we can also  use {\it inclusion-exclusion} (after applying Matsuo's bijection).
Procedure {\tt RIN(n,a,b)} uses {\it inclusion-exclusion} to compute,
more generally, in time $O(n^5)$, the first $n$ terms of the set of permutations $\pi$ such that $\pi_{i+1}-\pi_i \neq 1$
except possibly if $i=a$  or $p_i=b$, for {\it any} $1 \leq a,b \leq n-1$.

It is very likely that the three-variable function $RIN(n,a,b)$ is holonomic in all its variables
(i.e. satisfies linear recurrences with coefficients that are polynomials of $(n,a,b)$ in each of the three variables). 
A possible strategy would
be to guess it from empirical data and then prove that it satisfies some recurrences linking $RIN(n,a,b)$ to its neighbors,
obtained from combinatorial considerations, and then take the `diagonal' $RIN(n, \lfloor \frac{n+1}{2} \rfloor,  \lfloor \frac{n+1}{2} \rfloor )$,
that would have to be holonomic in $n$.

Using Matsuo's bijection we wrote the absolute value analog, that enabled us to generate $150$ terms for
the companion sequence $b_{2,2}(n)$, {\tt https://oeis.org/A110128}, see:

{\tt https://sites.math.rutgers.edu/\~{}zeilberg/tokhniot/oResPerms22VrmA.txt} \quad .

With some optimization, and a larger computer, this can probably be extended to at least $300$ terms.
Using an adaptation of his original approach, Rintaro Matsuo generated $600$ terms, and this enabled
Manuel Kauers to conjecture a linear recurrence of order $24$ and degree $40$ satisfied by this sequence,
giving get more evidence that the sequence $a_{r,s}(n)$ and $b_{r,s}(n)$ are always $P$-recursive.

{\eightrm [To see Kauers' conjecture based on Matsuo's data type {\tt Ope22Kv(n,N)} in the Maple package {\tt ResPermsRM.txt}]} \quad .

Matsuo's bijection, and its natural extensions, suggest a whole new research direction. Given finite {\bf sets} of positive integers $A$ and $B$,
efficiently compute, and study, the sequences, let's call them $R_{A,B}(n)$,
of permutations of $\{1, \dots , n\}$ that $\pi_{i+1}-\pi_i \neq 1$ except that  $\pi_{i+1}-\pi_i= 1$ is allowed if $i \in A$ or $\pi_i \in B$.
One can also make the elements of $A$ and $B$  depend on $n$.

Also of interest is the absolute value version, i.e.
counting  permutations of $\{1, \dots , n\}$ such that $|\pi_{i+1}-\pi_i| \neq 1$ {\it except} that  $|\pi_{i+1}-\pi_i|= 1$ is allowed, whenever $i \in A$ or $\pi_i \in B$,
for $n \geq max(A \cup B)$. One can use an extension of Matsuo's approach, or our inclusion-exclusion version (or both).
Either way should give a
{\it polynomial time} algorithms. Alas, the implied {\bf degree} in the `polynomial time' would be eventually
prohibitive. More importantly, this suggests, but by no means proves, that these sequences are indeed always $P$-recursive (holonomic), i.e. satisfy
a linear recurrence equation with polynomial coefficients.

We will be very happy if you, dear readers, will pursue all these ramifications.
In particular design and implement  $O(n^6)$ algorithms to compute $a_{3,3}(n)$ and $b_{3,3}(n)$.

{\bf References}

[AZ] Gert Almkvist and Doron Zeilberger, {\it The method of differentiating under the integral sign},
 J. Symbolic Computation {\bf 10} (1990), 571-591. \hfill\break
{\tt https://sites.math.rutgers.edu/\~{}zeilberg/mamarim/mamarimPDF/duis.pdf} \quad .

[KK] Manuel Kauers and Christoph Koutschan, {\it Guessing with little Data}, Proceedings of ISSAC'22, 83-90.
\quad {\tt https://arxiv.org/abs/2202.07966} \quad .

[N] Enrique Navarrete,  {\it Generalized K-shift forbidden substrings in permutations}, \hfill\break
{\tt https://arxiv.org/pdf/1610.06217.pdf}

[P] P. Poulet, {\it  Query 4750: Permutations}, L'Interm\'ediaire des Math\'ematiciens {\bf 26} (1919), 117-121.  \hfill\break
[See the OEIS entry {\tt https://oeis.org/A002464} for scans of its pages]

[Rio] John Riordan, {\it A recurrence for permutations without rising or falling successions},
Ann. Math. Statist. {\bf 36} (1965), 708-710.

[Rin] Rintaro Matsuo, {\it $O(n^4)$ code to calculate $a(n)$}, \hfill\break
[It is linked from {\tt https://oeis.org/A189281}] \quad .

[Ro] David P. Robbins, {\it The probability that neighbors remain neighbors after random rearrangements},
Amer. Math. Monthly {\bf 87} (1980), 122-124. \hfill\break
{\tt https://www.jstor.org/stable/2321990} \quad .

[S11] N. J. A. Sloane, {\it The On-Line Encyclopedia of Integer Sequences}, sequence A255,  \hfill\break
{\tt https://oeis.org/A000255}.

[S11a] N. J. A. Sloane, {\it The On-Line Encyclopedia of Integer Sequences}, sequence A2464,  \hfill\break
{\tt https://oeis.org/A002464}.

[S12] N. J. A. Sloane, {\it The On-Line Encyclopedia of Integer Sequences}, sequence A55790, \hfill\break
{\tt  https://oeis.org/A055790}.

[S13] N. J. A. Sloane, {\it The On-Line Encyclopedia of Integer Sequences}, sequence  A277609, \hfill\break
{\tt  https://oeis.org/A277609}.

[S14] N. J. A. Sloane, {\it The On-Line Encyclopedia of Integer Sequences}, sequence  A277563, \hfill\break
{\tt https://oeis.org/A277563 }.

[S15] N. J. A. Sloane, {\it The On-Line Encyclopedia of Integer Sequences}, sequence  A280425, \hfill\break
{\tt  https://oeis.org/A280425 }.

[S22] N. J. A. Sloane, {\it The On-Line Encyclopedia of Integer Sequences}, sequence A189281, \hfill\break
{\tt https://oeis.org/A189281}.

[S22a] N. J. A. Sloane, {\it The On-Line Encyclopedia of Integer Sequences}, sequence A110128, \hfill\break
{\tt https://oeis.org/A110128}.

[S33] N. J. A. Sloane, {\it The On-Line Encyclopedia of Integer Sequences}, sequence A189282, \hfill\break
{\tt https://oeis.org/A189282}.

[S33a] N. J. A. Sloane, {\it The On-Line Encyclopedia of Integer Sequences}, sequence A117574, \hfill\break
{\tt https://oeis.org/A117574}.

[S44] N. J. A. Sloane, {\it The On-Line Encyclopedia of Integer Sequences}, sequence A189283, \hfill\break
{\tt  https://oeis.org/A189283}.

[S44a] N. J. A. Sloane, {\it The On-Line Encyclopedia of Integer Sequences}, sequence A189255, \hfill\break
https://oeis.org/A189255

[S55] N. J. A. Sloane, {\it The On-Line Encyclopedia of Integer Sequences}, sequence A189284, \hfill\break
{\tt  https://oeis.org/A189284}.

[S55a] N. J. A. Sloane, {\it The On-Line Encyclopedia of Integer Sequences}, sequence A189256, \hfill\break
{\tt  https://oeis.org/A189256}.

[S66] N. J. A. Sloane, {\it The On-Line Encyclopedia of Integer Sequences}, sequence A189285, \hfill\break
{\tt  https://oeis.org/A189285}.

[S66a] N. J. A. Sloane, {\it The On-Line Encyclopedia of Integer Sequences}, sequence A189271, \hfill\break
{\tt  https://oeis.org/A189271}.

[T] Roberto Tauraso, {\it The dinner table problem: the rectangular Case}, INTEGERS, {\bf 6} (2006), paper A11. \quad
{\tt https://www.emis.de/journals/INTEGERS/papers/g11/g11.pdf} \quad.

[W] Herbert S. Wilf, {\it What is an answer?},  American Mathematical Monthly
{\bf 89} (1982),  289-292.

[Z1] Doron Zeilberger, {\it The method of creative telescoping}, J. Symbolic Computation {\bf 11}(1991) 195-204, \quad
{\tt http://sites.math.rutgers.edu/\~{}zeilberg/mamarimY/CreativeTelescoping.pdf}

[Z2]  Doron Zeilberger, {\it Automatic CounTilings}, Personal Journal of Shalosh B. Ekhad and Doron Zeilberger, Jan. 20, 2006. \hfill\break
{\tt https://sites.math.rutgers.edu/\~{}zeilberg/mamarim/mamarimhtml/tilings.html} \quad .

[Z3] Doron Zeilberger, {\it Automatic enumeration of generalized M\'enage numbers}, S\'eminaire Lotharingien de Combinatoire, {\bf 71} (2014), article B71a.  \hfill\break
{\tt https://sites.math.rutgers.edu/\~{}zeilberg/mamarim/mamarimhtml/menages.html} \quad .

\bigskip
\hrule
\bigskip
George Spahn and Doron Zeilberger, Department of Mathematics, Rutgers University (New Brunswick), Hill Center-Busch Campus, 110 Frelinghuysen
Rd., Piscataway, NJ 08854-8019, USA. \hfill\break
Email: {\tt  gs828 at math dot rutgers dot edu} \quad, \quad {\tt DoronZeil] at gmail dot com}   \quad .

First Written: {\bf Nov. 3, 2022}. This version: {\bf Dec. 5, 2022}.

\end